\documentclass[a4paper,12pt]{amsart}
\usepackage[latin1]{inputenc}
\usepackage[francais]{babel}
\usepackage[T1]{fontenc}
\usepackage{latexsym}
\usepackage{amsmath,amsthm}
\usepackage{amsfonts}
\usepackage{epsfig}
\usepackage{amssymb}
\usepackage{amscd}
\usepackage[all]{xy}
\usepackage{stmaryrd}
\input xy

\newtheorem{Lemma}{Lemme}
\newtheorem{Theorem}[Lemma]{Th{\'e}or{\`e}me}

\newtheorem{Proposition}[Lemma]{Proposition}

\theoremstyle{definition}
\newtheorem{Def}[Lemma]{D{\'e}finition}

\theoremstyle{remark}
\newtheorem{Remark}[Lemma]{Remarque}

\newcommand{\RR}{\mathbb{R}}

\newcommand{\A}{\mathcal{A}}

\newcommand{\abs}[1]{\lvert*1\rvert}

\newcommand{\blankbox}[2]{%
  \parbox{\columnwidth}{\centering
    \setlength{\fboxsep}{0pt}%
    \fbox{\raisebox{0pt}[*2]{\hspace{*1}}}%
  }%
}

\title[Propri{\'e}t{\'e} (T) tordue]{Sur la Propri{\'e}t{\'e} (T) tordue par un produit tensoriel}

\author{Maria-Paula Gomez-Aparicio}

\address{Maria-Paula {\sc Gomez-Aparicio} : Institut de Math{\'e}matiques de Jussieu, Equipe d'alg{\`e}bres d'Op{\'e}rateurs,
175 rue du chevaleret, 75013 Paris, France.}

\email{gomez@math.jussieu.fr}

\date{}

\keywords{Unitary representation, matrix coefficients, K-types}

\begin{document}

\maketitle
\begin{abstract}
Dans cet article on d{\'e}finit une version tordue de la propri{\'e}t{\'e} (T) de Kazhdan en consid{\'e}rant des produits
tensoriels de repr{\'e}sentations unitaires par des repr{\'e}sentations
irr{\'e}ductibles de dimension finie non-unitaires d'un groupe topologique.
On d{\'e}montre, en utilisant la d{\'e}croissance uniforme des coefficients de matrices des repr{\'e}sentations unitaires d'un groupe de Lie
r{\'e}el simple $G$, ayant la propri{\'e}t{\'e} (T) de Kazhdan,  que toute repr{\'e}sentation irr{\'e}ductible de dimension finie $\rho$ de $G$
est isol{\'e}e parmi les repr{\'e}sentations de la forme $\rho\otimes\pi$, o{\`u}
$\pi$  parcourt les repr{\'e}sentations unitaires irr{\'e}ductibles de $G$, dans un
sens que l'on précisera.

\vspace{2mm}
\noindent{\sc {\large a}bstract.} In this article, we consider tensor
products of unitary representations by irreducible non-unitary finite
dimensional representations of topological groups to define a property
that is a twisting of Kazhdan's Property (T). We use the uniform
decay of the matrix coefficients of unitary representations, to show
that for most of the real semi-simple Lie groups having Kazhdan's
Property (T), any finite dimensional
irreducible representation  $\rho$ of $G$, is isolated among representations of the form $\rho\otimes\pi$, where $\pi$ ranges over the irreducible unitary representations, in a sense to be made precise.
\end{abstract}

\section*{Introduction}
Un groupe topologique $G$ a la propri{\'e}t{\'e} (T) de Kazhdan si sa
repr{\'e}sentation triviale est isol{\'e}e dans son dual unitaire,
$\widehat{G}$. En 1967, Kazhdan a introduit cette propri{\'e}t{\'e} pour
{\'e}tudier la structure des r{\'e}seaux dans les groupes de Lie r{\'e}els.
Il d{\'e}montre que tout r{\'e}seau d'un groupe de Lie r{\'e}el $G$ qui a la
propri{\'e}t{\'e} (T), est de type fini \cite{Kazhdan}. Plus tard, en 1981, Akemann et Walter on donné une caractérisation $C^*$-algébrique de la propriété (T) \cite[Lemma 2]{Akemann-Walter}. Ils démontrent qu'un groupe topologique localement compact $G$ a la propri{\'e}t{\'e}
(T) si et seulement s'il existe un idempotent $p$ dans la
$C^*$-alg{\`e}bre maximale $C^*(G)$ de $G$ tel que pour toute repr{\'e}sentation $(\pi,H)$
unitaire de $G$, $\pi(p)$ soit la projection orthogonale sur le sous-espace de
$H$ form{\'e} des
vecteurs $G$-invariants (voir aussi \cite[Proposition 2]{Valette}).\\
Dans cet article, on va
d{\'e}finir un renforcement de la propri{\'e}t{\'e} (T) en consid{\'e}rant des produits
tensoriels de repr{\'e}sentations unitaires par des repr{\'e}sentations irr{\'e}ductibles non unitaires de
dimension finie. Pour ceci, on va définir, pour toute représentation irréductible de dimension finie d'un groupe topologique, un analogue \emph{tordu} de la $C^*$-algèbre maximale de $G$, que l'on va noter $\mathcal{A}_G$, et on va définir une propriété (T) \emph{tordue} en termes de celle-ci.
\begin{Def}
 Soient $G$ un groupe localement compact et $(\rho,V)$ une repr{\'e}sentation irr{\'e}ductible de
 dimension finie de $G$. Soit $C_{c}(G)$ l'espace vectoriel des fonctions continues {\`a} support
 compact sur $G$ et soit $\mathcal{A}_G$ la compl{\'e}tion de
$C_{c}(G)$ pour la norme $\|. \|_{\mathcal{A}_G}$ donn{\'e}e par:
$$\|f\|_{\mathcal{A}_G}=\sup\limits_{(\pi,H_{\pi})}\|(\rho\otimes\pi)(f)\|_{\mathcal{L}(V\otimes
H_{\pi})},$$ pour $f\in C_{c}(G)$ et o{\`u} $(\pi,H_{\pi})$ est une repr{\'e}sentation unitaire de $G$.\\
On dira que $G$
a la propri{\'e}t{\'e} (T) \emph{tordue par $\rho$} s'il existe un idempotent
$p_G$ dans $\A_G$ tel que: $\rho(p_G)=\mathrm{Id}_{V}$
et, pour toute $\pi$ repr{\'e}sentation unitaire de $G$ qui ne
contient pas la repr{\'e}sentation triviale, $(\rho\otimes\pi)(p_G)=0$.
\end{Def}

Si un groupe $G$ a la propri{\'e}t{\'e} (T) tordue par $\rho$,
on dira alors que $\rho$ est \emph{isol{\'e}e parmi les repr{\'e}sentation de la forme $\rho\otimes\pi$}, o{\`u} $\pi$ est
une repr{\'e}sentation unitaire irr{\'e}ductible de $G$, terminologie que l'on
justifiera (voir la proposition \ref{T}).\\

De plus, pour toute représentation irréductible de dimension finie $\rho$ on va d{\'e}finir, de la
m{\^e}me fa{\c c}on, un analogue \emph{tordu}
de $C^*_r(G)$, la $C^*$-alg{\`e}bre r{\'e}duite de $G$, que l'on notera
$\mathcal{A}_G^r$, et on d{\'e}montrera que
si $G$ a la propri{\'e}t{\'e} (T) tordue $\rho$, alors les alg{\`e}bres $\mathcal{A}_G$ et $\mathcal{A}_G^r$ n'ont pas la m{\^e}me
K-th{\'e}orie, l'int{\'e}r{\^e}t pour nous {\'e}tant de calculer la K-th{\'e}orie de ces alg{\`e}bres.\\
On d{\'e}montre aussi que la propri{\'e}t{\'e} (T) tordue par une repr{\'e}sentation
$\rho$ est h{\'e}rit{\'e}e par tout r{\'e}seau cocompact du groupe: Si $G$ a (T) tordue par
$\rho$ et $\Gamma$ est un r{\'e}seau cocompact de $G$, alors $\Gamma$ a la propri{\'e}t{\'e} (T) tordue par
$\rho|_{\Gamma}$.\\
On ne sait pas encore si tout groupe localement compact ayant la propriété (T) a aussi la propriété (T) tordue par toute représentation irréductible de dimension finie, mais la derni{\`e}re partie de l'article est consacr{\'e}e {\`a} la d{\'e}monstration du
fait que, au moins, beaucoup de groupes de Lie ayant la propri{\'e}t{\'e} (T) ont aussi
la propri{\'e}t{\'e} (T) tordue par n'importe quelle repr{\'e}sentation
irr{\'e}ductible de dimension finie. On sait que tout groupe de Lie $G$ r{\'e}el simple connexe de centre fini de rang r{\'e}el $\geq 2$
ou localement isomorphe {\`a} $Sp(n,1)$ pour $n\geq 2$ ou {\`a} $F_{4(-20)}$, a la
propri{\'e}t{\'e} (T) \cite{de la Harpe-Valette}. Plus fort encore, il v{\'e}rifie une d{\'e}croissance
uniforme des coefficients
de matrice des repr{\'e}sentations unitaires de $G$ qui n'ont pas de vecteurs invariants non nuls \cite{Cowling}.
On utilise cette propri{\'e}t{\'e} donn{\'e}e par le th{\'e}or{\`e}me \ref{Howe}, pour
montrer le résultat suivant:
\begin{Theorem}\label{intro}
Soit $G$ un groupe de Lie r{\'e}el alg{\'e}brique
semi-simple de centre fini, connexe et simplement connexe (au sens
alg{\'e}brique), et tel que chaque facteur simple de $G$ est ou bien de
rang r{\'e}el sup{\'e}rieur ou {\'e}gal {\`a} $2$, ou bien localement isomorphe {\`a}
$Sp(n,1)$ pour $n\geq 2$ ou
{\`a} $F_{4(-20)}$, et soit $\rho$ une repr{\'e}sentation
irr{\'e}ductible de dimension finie de $G$. Alors $G$ a la propri{\'e}t{\'e} (T) tordue par $\rho$.

\end{Theorem}

Le cas o{\`u} $G=SL_m(\mathbb{R})$ et $m\geq 3$ a {\'e}t{\'e} pr{\'e}sent{\'e} par l'auteur lors de l'{\'e}cole d'{\'e}t{\'e} intitul{\'e}e "Topological and Geometric Methods for
Quantum Field Theory" pendant l'{\'e}t{\'e} 2005 {\`a} Villa de Leyva en Colombie et va {\^e}tre publi{\'e} dans les comptes-rendus de celle-ci \cite{Gomez}.\\

David Fisher et Theron Hitchman d{\'e}finissent dans \cite{Fisher} une
propri{\'e}t{\'e}, qu'ils appellent $F\otimes H$, en termes de $1$-cohomologie, qui ressemble {\`a} la propri{\'e}t{\'e} (T) tordue par $\rho$, mais il
n'y a pas d'implication entre les deux propri{\'e}t{\'e}s.\\

\medskip
\noindent {\bf Remerciements.} Je voudrais remercier Vincent
Lafforgue pour ses nombreuses suggestions et sa grande
disponibilit{\'e} et Bachir Bekka pour ses {\'e}claircissements dans le
cas des groupes de rang 1.

\section{Propri{\'e}t{\'e} (T) tordue}
\noindent{\bf D{\'e}finitions et terminologie.}
On rappelle la d{\'e}finition de propri{\'e}t{\'e} (T) \emph{tordue} énonc{\'e}e dans l'introduction:
\begin{Def}\label{tordue}
 Soient $G$ un groupe localement compact et $(\rho,V)$ une repr{\'e}sentation irr{\'e}ductible de
 dimension finie de $G$, o{\`u} $V$ est un espace vectoriel complexe muni d'une norme hermitienne. Soit $C_{c}(G)$ l'espace vectoriel des fonctions continues {\`a} support
 compact sur $G$ et soit $\mathcal{A}_G$ la compl{\'e}tion de
$C_{c}(G)$ pour la norme $\|. \|_{\mathcal{A}_G}$ donn{\'e}e par:
$$\|f\|_{\mathcal{A}_G}=\sup\limits_{(\pi,H_{\pi})}\|(\rho\otimes\pi)(f)\|_{\mathcal{L}(V\otimes
H_{\pi})},$$ pour $f\in C_{c}(G)$, o{\`u} $(\pi,H_{\pi})$ varie parmi
les repr{\'e}sentations unitaires de $G$.\\
On note $1_G$ la repr{\'e}sentation triviale de $G$ et on dit que $G$
a la propri{\'e}t{\'e} (T) tordue par $\rho$ s'il existe un idempotent
$p_G$ dans $\A_G$ tel que:\\
$\rho(p_G)=\mathrm{Id}_{V}$
et, pour toute $\pi$ repr{\'e}sentation unitaire de $G$ qui ne
contient pas $1_G$, $(\rho\otimes\pi)(p_G)=0$.
\end{Def}

\begin{Remark}
L'alg{\`e}bre $\mathcal{A}_G$ est une alg{\`e}bre de Banach involutive et toute repr{\'e}sentation du groupe $G$ de la forme $\rho\otimes\pi$ avec $\pi$ unitaire peut {\^e}tre prolong{\'e}e,
 de fa{\c c}on {\'e}vidente, en une repr{\'e}sentation de $\A_G$ que l'on note aussi, par abus de
 notation, $\rho\otimes\pi$.
\end{Remark}

\begin{Remark}
 Si $G$ a la propri{\'e}t{\'e} (T) tordue par $(\rho,V)$, alors, pour toute $(\pi,H)$ repr{\'e}sentation unitaire de $G$,
 $(\rho\otimes\pi)(p_G)$ est la projection de $V\otimes H$ sur
 $V\otimes H^G$ parallèlement {\`a} $V\otimes (H^G)^{\bot}$, o{\`u} $H^G$ est le sous-espace de $H$
 form{\'e} des vecteurs $G$-invariants. En effet, il suffit d'{\'e}crire $\pi$
 de la forme
 $\pi_0\oplus\pi_1$, o{\`u} $\pi_1$ est la sous-repr{\'e}sentation de $\pi$ sur
 $(H^G)^{\bot}$ qui ne contient
 pas la repr{\'e}sentation triviale et $\pi_0$ est la
 sous-repr{\'e}sentation de $\pi$ qui a pour espace $H^{G}$ et qui est
 {\'e}quivalente {\`a} $1_G$.
\end{Remark}

\begin{Def}
 Soit $G$ un groupe localement compact, $\widehat{G}$ son dual unitaire et $(\rho,V)$
 une repr{\'e}sentation irr{\'e}ductible de $G$. On d{\'e}finit une nouvelle topologie sur
 $\widehat{G}$, que l'on appelle \emph{tordue par $\rho$},
  de la mani{\`e}re suivante: Si $P$ est un sous-ensemble de $\widehat{G}$ et
 $\pi\in\widehat{G}$
  alors $\pi$
 est dans l'adh{\'e}rence tordue par $\rho$ de $P$, si
  $\rho\otimes\pi$ est contenu dans l'adh{\'e}rence de Fell de
 $\rho\otimes P:=\{\rho\otimes\pi'| \pi'\in P\}$ dans le dual de $\A_G$ ( voir \cite{Fell} Chapitre VII ).\\
On note $\widehat{G}^{\rho}$ l'espace $\widehat{G}$ muni de cette
topologie.
\end{Def}

\begin{Proposition}\label{T}
 Soit $G$ un groupe topologique localement compact et $\rho$ une repr{\'e}sentation
 irr{\'e}ductible de dimension finie de $G$. Si $G$ a la propri{\'e}t{\'e} (T) tordue par $\rho$ alors la
 repr{\'e}sentation triviale de $G$ est isol{\'e}e dans $\widehat{G}^{\rho}$.
\end{Proposition}

\begin{proof}
 \sloppy Supposons que $G$ a la propri{\'e}t{\'e} (T) tordue par $\rho$ et soit
 $p_G$ l'idempotent dans $\mathcal{A}_G$ qui v{\'e}rifie les
 propri{\'e}t{\'e}s donn{\'e}es par la d{\'e}finition \ref{tordue}. Supposons que $\rho$ soit contenue dans l'adh{\'e}rence de Fell de l'ensemble $\{\rho\otimes\pi | \pi\,\hbox{est une repr{\'e}sentation unitaire}\}$ dans le dual de $\A_G$.
 On a alors
 $$\bigcap_{1_G\nsubseteq\pi}\mathrm{Ker}(\rho\otimes\pi)\subset\mathrm{Ker}(\rho),$$
(voir par exemple \cite[Chapitre VII Proposition 3.11]{Fell}).\\
 Or,
 $p_G\in\bigcap_{1_G\nsubseteq\pi}\mathrm{Ker}(\rho\otimes\pi)$ et $\rho(p_G)\neq 0$, et
 donc $\rho$ est isol{\'e}e, pour la topologie de Fell, dans le sous-ensemble du dual de $\A_G$
 form{\'e} des repr{\'e}sentations de la forme $\rho\otimes\pi$, o{\`u} $\pi$ est une repr{\'e}sentation unitaire de $G$.
Par cons{\'e}quent, $1_G$ est isol{\'e}e dans $\widehat{G}^{\rho}$.
\end{proof}
Cette proposition justifie la terminologie utilis{\'e}e:\\
Si un groupe topologique $G$ a la propri{\'e}t{\'e} (T) tordue par une repr{\'e}sentation irr{\'e}ductible de dimension finie $\rho$, on dira alors que $\rho$ est \emph{isol{\'e}e parmi les repr{\'e}sentations de la forme $\rho\otimes\pi$}, o{\`u} $\pi$ est  une repr{\'e}sentation unitaire irr{\'e}ductible de $G$.
\begin{Remark}
Si $\rho=1_G$, alors $G$ a la propri{\'e}t{\'e} (T) tordue par $1_G$ s'il
existe un idempotent $p_G$ dans la $C^*$-alg{\`e}bre maximale de $G$,
tel que pour toute repr{\'e}sentation unitaire $\pi$ de $G$,
$\pi(p_G)$ soit la projection sur l'espace des vecteurs
$G$-invariants. Dans ce cas, la proposition \ref{T} est une
{\'e}quivalence et c'est le r{\'e}sultat connu qui dit qu'un groupe
localement compact a la propri{\'e}t{\'e} (T) usuelle si et seulement si $C^*(G)$ s'{\'e}crit
comme une somme directe de $C^*$-alg{\`e}bres de la forme:
$$C^*(G)=\mathrm{Ker}(1_G)\oplus I,$$o{\`u} $I$ est un id{\'e}al bilat{\`e}re
ferm{\'e} de $C^*(G)$ \cite{Akemann-Walter}, \cite{Valette}.
\end{Remark}

L'objectif principal de cet article est de prouver que beaucoup de
groupes de Lie qui ont la propri{\'e}t{\'e} (T), v{\'e}rifient aussi la
propri{\'e}t{\'e} (T) tordue par n'importe quelle repr{\'e}sentation irr{\'e}ductible de
dimension finie.\\
En
s'inspirant de la remarque pr{\'e}c{\'e}dente, pour tout groupe topologique localement compact $G$ et toute repr{\'e}sentation irr{\'e}ductible de dimension finie $\rho$ de $G$, on consid{\`e}re l'espace
vectoriel $C_c(G)$ form{\'e} des fonctions continues {\`a} support compact
sur $G$ et on d{\'e}finit deux nouvelles compl{\'e}tions de celui-ci de la mani{\`e}re suivante:\\

Soit $\mathcal{A'}_G$ la compl{\'e}tion de $C_{c}(G)$ par rapport {\`a} la
norme $\|. \|_{\mathcal{A'}_G}$ donn{\'e}e par:
$$\|f\|_{\mathcal{A'}_G}=\sup\limits_{(\pi,H_{\pi})}\|(\rho\otimes\pi)(f)\|_{\mathcal{L}(V\otimes
H_{\pi})},$$ pour $f\in C_{c}(G)$, o{\`u} $(\pi,H_{\pi})$ varie parmi les repr{\'e}sentations unitaires de $G$ qui ne contiennent pas la repr{\'e}sentation triviale.\\
Et soit $\mathcal{A}''_G$ la compl{\'e}tion de $C_{c}(G)$
pour la norme $\|. \|_{\mathcal{A''}_G}$ donn{\'e}e par:
$$\|f\|_{\mathcal{A}''_G}=\|\rho(f)\|_{\mathrm{End}(V)},$$
pour tout $f\in C_{c}(G)$.\\
On remarque alors qu'on a deux morphismes d'alg{\`e}bres de Banach:
$$\begin{array}{lcl}
\Theta_{1}:\mathcal{A}_G\rightarrow\mathcal{A'}_G &\!\!\!\mathrm{and}\!\!\!& \Theta_{2}:\mathcal{A}_G\rightarrow\mathcal{A''}_G\\
\end{array}.$$
Soit
$\Theta:\mathcal{A}_G\rightarrow\mathcal{A'}_G\oplus\mathcal{A''}_G$ le
prolongement {\`a} $\mathcal{A}_G$ du morphisme
donn{\'e} sur $C_c(G)$ par: $\Theta(f)=(\Theta_{1}(f),\Theta_{2}(f))$. C'est un morphisme d'alg{\`e}bres de Banach.\\

On a alors l'{\'e}quivalence {\'e}vidente suivante:
\begin{Proposition}\label{evidente}
Le groupe $G$ a la propri{\'e}t{\'e} (T)
tordue par $\rho$ si et seulement si le morphisme d'alg{\`e}bres de
Banach $\Theta$ est un isomorphisme.\\
\end{Proposition}

\noindent{\bf Relation avec la K-th{\'e}orie}
Soient $G$ un groupe topologique et $(\rho,V)$ une repr{\'e}sentation
irr{\'e}ductible de dimension finie de $G$. De la m{\^e}me fa{\c c}on que l'on a
d{\'e}fini l'alg{\`e}bre $\mathcal{A}_G$, qui est l'analogue \emph{tordu par $\rho$} de la
$C^*$-alg{\`e}bre maximale de $G$, on peut d{\'e}finir une alg{\`e}bre de Banach
\emph{r{\'e}duite tordue par $\rho$} de $G$ comme étant la compl{\'e}tion de $C_c(G)$
pour la norme:
$$\|f\|_{\mathcal{A}_G^r}=\|\rho\otimes L_G(f)\|_{\mathcal{L}(V\otimes\mathrm{L}^2(G))},$$
pour $f\in C_c(G)$ et o{\`u} $ L_G$ est la repr{\'e}sentation r{\'e}guli{\`e}re
gauche de
$G$ dans $\mathrm{L}^2(G)$.  On note $\mathcal{A}_G^r$ cette compl{\'e}tion.\\
On a alors un unique morphisme d'alg{\`e}bres de Banach prolongeant l'identit{\'e}
sur $C_c(G)$:
$$\rho\otimes L_G:\mathcal{A}_G\rightarrow\mathcal{A}_G^r,$$
qui d{\'e}finit un morphisme en K-th{\'e}orie
$$(\rho\otimes L_G)^*:\mathrm{K}(\mathcal{A}_G)\rightarrow
\mathrm{K}(\mathcal{A}_G^r).$$
\begin{Proposition}
Si $G$ a la propri{\'e}t{\'e} (T) tordue par $\rho$ et $G$ n'est pas un groupe compact
 alors les alg{\`e}bres
$\mathcal{A}_G$ et $\mathcal{A}_G^r$ n'ont pas la m{\^e}me K-th{\'e}orie,
c'est-{\`a}-dire que $(\rho\otimes L_G)^*$
n'est pas un isomorphisme.
 \end{Proposition}
\begin{proof}
Supposons que $G$ soit un groupe localement compact non-compact ayant la propri{\'e}t{\'e} (T)
tordue par une représentation $\rho$. Il existe alors un idempotent non-nul $p_G\in\mathcal{A}_G$
tel que $\rho(p_G)=\mathrm{Id}_{V}$ et, pour toute repr{\'e}sentation $\pi$ unitaire de $G$ qui ne contient pas
la repr{\'e}sentation triviale , $(\rho\otimes\pi)(p_G)=0$. 
Comme $G$ n'est pas compact on peut prendre $\pi=L_G$ et, par cons{\'e}quent,
$(\rho\otimes L_G)(p_G)=0$, ce qui montre que
$(\rho\otimes L_G)^*$ n'est pas un morphisme injectif.\\
\end{proof}

\noindent{\bf Propri{\'e}t{\'e} d'h{\'e}r{\'e}dit{\'e}} On va démontrer que la propriété (T) tordue par une représentation irréductible de dimension finie est héritée par les réseaux cocompacts du groupe. Pour qu'un énoncé de ce style ait un sens, on a d'abord besoin du lemme suivant:
\begin{Lemma}
Si $\Gamma$ est un r{\'e}seau cocompact d'un groupe topologique localement compact $G$ ayant la
propri{\'e}t{\'e} (T) tordue par une repr{\'e}sentation irr{\'e}ductible de dimension
finie $(\rho, V)$, alors
$\rho|_{\Gamma}$ est une repr{\'e}sentation irr{\'e}ductible de $\Gamma$.
\end{Lemma}

\begin{proof}
\sloppy Soit $\mathrm{Hom}_{\Gamma}(V|_{\Gamma},V|_{\Gamma})$ l'ensemble des
morphismes
$\Gamma$-invariants de $V|_{\Gamma}$ dans $V|_{\Gamma}$. On a un
morphisme injectif de
$\mathrm{Hom}_{\Gamma}(V|_{\Gamma},V|_{\Gamma})$ dans
$\mathrm{Hom}_G(V,V\otimes \mathrm{L^2}(G/\Gamma))$, o{\`u}
$L^2(G/\Gamma)$ est l'espace de la repr{\'e}sentation r{\'e}guli{\`e}re
quasi-invariante de $G$. En effet, soit
$T\in\mathrm{Hom}_{\Gamma}(V|_{\Gamma},V|_{\Gamma})$ et $F_T$ la
fonction continue sur $G$ {\`a} valeurs dans $\mathrm{End}(V)$ telle que
$F_T(g)=\rho(g)T\rho(g)^{-1}$, pour tout $g\in G$. Comme $T$ est
$\Gamma$-{\'e}quivariant, $F_T$ est une fonction continue (qui est en plus $G$-{\'e}quivariante) sur $G/\Gamma$ {\`a}
valeurs dans $\mathrm{End}(V)$. Comme
$G/\Gamma$ est compact, $F_T$ appartient {\`a}
$\mathrm{L^2}(G/\Gamma,\mathrm{End}(V))$.\\
Soit $$F(T):V\rightarrow
\mathrm{L^2}(G/\Gamma,V),$$ tel que
$F(T)(v)(x)=F_T(x)v$, pour tout $v\in V$ et tout $x\in G/\Gamma$. C'est une application
lin{\'e}aire continue $G$-{\'e}quivariante, donc
$F(T)\in\mathrm{Hom}_G(V,\mathrm{L^2}(G/\Gamma,V))=\mathrm{Hom}_G(V,\mathrm{L^2}(G/\Gamma)\otimes
V)$.\\
De plus, si $1\in G$ est l'identit{\'e} de $G$, alors  $F_T(1)=T$ et donc, la
correspondance $T\mapsto F(T)$ d{\'e}finit un morphisme injectif de
$\mathrm{Hom}_{\Gamma}(V|_{\Gamma},V|_{\Gamma})$ dans
$\mathrm{Hom}_G(V,V\otimes \mathrm{L^2}(G/\Gamma))$.\\
Ceci implique que,
 $$\mathrm{dim}_{\mathbb{C}}(\mathrm{Hom}_{\Gamma}(V|_{\Gamma},V|_{\Gamma}))\leq\mathrm{dim}_{\mathbb{C}}(\mathrm{Hom}_G(V,V\otimes \mathrm{L^2}(G/\Gamma))).$$
Mais le fait que le groupe $G$ ait la propri{\'e}t{\'e} (T) tordue par $\rho$ implique que
$$\mathrm{Hom}_G(V,V\otimes
\mathrm{L^2}(G/\Gamma))=\mathrm{Hom}_G(V,V),$$
car, en effet, on peut {\'e}crire
$$\mathrm{L^2}(G/\Gamma)=\mathrm{L^2}(G/\Gamma)_0\oplus\mathrm{L^2}(G/\Gamma)_1,$$
o{\`u} $\mathrm{L^2}(G/\Gamma)_0$ est la sous-repr{\'e}sentation triviale
de $L^2(G/\Gamma)$ engendr{\'e}e par la fonction constante {\'e}gale {\`a} $1$ et $\mathrm{L^2}(G/\Gamma)_1$
est son orthogonal, qui ne contient pas $1_G$. Si en plus $G$ a (T) tordue par $\rho$,
$$\mathrm{Hom}_G(V,V\otimes\mathrm{L^2}(G/\Gamma)_1)=0,$$
d'o{\`u},
\begin{align*}
\mathrm{Hom}_G(V,V\otimes
\mathrm{L^2}(G/\Gamma))&=\mathrm{Hom}_G(V,V)\\
&=\mathbb{C}.\mathrm{Id}_{V},
\end{align*}
car $V$ est une
repr{\'e}sentation irr{\'e}ductible de $G$, d'o{\`u} le lemme.

\end{proof}
\begin{Proposition}
Soit $G$ un groupe localement compact et $\Gamma$ un r{\'e}seau
cocompact de $G$. Soit $\rho$ une repr{\'e}sentation irr{\'e}ductible de
dimension finie de $G$. Si $G$ a la propri{\'e}t{\'e} (T) tordue par
$\rho$ alors $\Gamma$ a la propri{\'e}t{\'e} (T) tordue par
$\rho|_\Gamma$.

\end{Proposition}

\begin{proof}
Supposons que $G$ a la propri{\'e}t{\'e} (T) tordue par $\rho$ et soit
$\Gamma$ un r{\'e}seau cocompact de $G$. L'existence de $\Gamma$ implique
que $G$ est unimodulaire. Soit $dg$ une mesure de Haar sur $G$ telle
que $G/\Gamma$ soit de mesure $1$. Par abus de notation,
 on notera de la m{\^e}me fa{\c c}on la repr{\'e}sentation $\rho$ de $G$ et sa restriction {\`a} $\Gamma$.
On rappelle que l'on note
$\mathcal{A}_G$ (resp. $\mathcal{A}_{\Gamma}$) la compl{\'e}tion de
$C_c(G)$ (resp. $C_c(\Gamma)$) pour la norme donn{\'e}e,
pour $f\in C_c(G)$ (resp. $f\in C_c(\Gamma)$), par:
$$\|f\|=\sup\limits_{(\pi,H_{\pi})}\|(\rho\otimes\pi)(f)\|_{\mathcal{L}(V\otimes
H_{\pi})},$$ o{\`u} le supremum est pris parmi les repr{\'e}sentations
unitaires de $G$
(resp. de $\Gamma$).\\
Si $(\pi,H)$ est une repr{\'e}sentation de $G$ (resp. de $\Gamma$) on
note $H^{G}$
 (resp. $H^{\Gamma}$) le sous-espace de $H$ form{\'e} des vecteurs invariants.\\
Soit $p_G$ l'idempotent de $\mathcal{A}_G$ tel que pour toute
$(\pi,H_{\pi})$
 repr{\'e}sentation unitaire de $G$, $(\rho\otimes\pi)(p_G)$ est la projection de
 $V\otimes H_{\pi}$ sur $V\otimes H_{\pi}^{G}$.\\
Supposons qu'il existe une fonction $f\in C_c(G)$ telle que
\begin{equation}\label{partition}
\sum_{\gamma\in\Gamma}f(g\gamma)\rho(\gamma)=\rho(g^{-1}),
\end{equation}
pour tout $g\in G$ (ce qui implique que
$\int_Gf(g)\rho(g)dg=\mathrm{Id}_{V}$).\\

 On veut construire un idempotent $p_{\Gamma}\in\mathcal{A}_{\Gamma}$,
 tel que pour toute repr{\'e}sentation $(\pi,H)$ unitaire
  de $\Gamma$, $(\rho\otimes\pi)(p_{\Gamma})$ soit la projection de $V\otimes H$ sur $V\otimes H^{\Gamma}$.\\
Soit une suite $(p_{G}^n)_n$ dans $C_c(G)$ qui converge vers $p_G$
dans $\mathcal{A}_G$, soit $(\pi,H)$ une repr{\'e}sentation unitaire de
$\Gamma$ et soit $(\mathrm{Ind}(\pi),\mathrm{Ind}(H))$ la
repr{\'e}sentation unitaire de $G$ obtenue par induction unitaire de
$\pi$. On note $\pi':=\mathrm{Ind}(\pi)$ et $H':=\mathrm{Ind}(H)$ pour
simplifier les notations et on remarque que $V\otimes H'$ est le compl{\'e}t{\'e}
de l'ensemble des fonction continues $s:G\rightarrow V\otimes H$ telles
que
$s(g\gamma)=(\mathrm{Id}_{\scriptstyle\mathrm{End}(V)}\otimes\pi')(\gamma^{-1})s(g)$,
pour $\gamma\in\Gamma$ et pour une norme $L^2$ sur $G/\Gamma$.\\
Soient $\alpha$ et $\beta$ les application lin{\'e}aires continues
suivantes:

\begin{align*}
 \alpha: V\otimes H &\rightarrow V\otimes H'\\
v\otimes\xi&\mapsto (g\mapsto \sum\limits_{\gamma\in\Gamma}\rho(g\gamma)v\otimes f(g\gamma)\pi(\gamma)\xi),\\
 \beta: V\otimes H' &\rightarrow V\otimes H\\
v\otimes s&\mapsto \int_{G}\rho(g^{-1})v\otimes f(g^{-1})s(g) dg.\\
\end{align*}

On v{\'e}rifie facilement que
$\beta\circ\alpha=\mathrm{Id}_{V\otimes H}$
et on a, pour tout $v\in V$ et tout $\xi\in H$,

\begin{align*}
\beta\circ(\rho\otimes\pi')(p_{G}^n)\circ&\alpha(v\otimes\xi)=\\
&\sum\limits_{\gamma\in\Gamma}\int_{G\times
G}p_G^n(g)f(x^{-1}g^{-1})f(x\gamma)(\rho\otimes\pi)(\gamma)(v\otimes\xi)dx
dg.
\end{align*}

Posons, pour $\gamma\in\Gamma$ et $n\in\mathbb{N}$,
 $$p_{\Gamma}^n(\gamma)=\int_{G\times G} p_G^n(g)f(x^{-1}g^{-1})f(x\gamma)dxdg.$$
La suite $(p_{\Gamma}^n)_n$ appartient {\`a} $C_c(\Gamma)$ et elle
converge dans $\mathcal{A}_\Gamma$ car:

$$\|p^n_{\Gamma}\|_{\mathcal{A}_{\Gamma}}=\sup\limits_{(\pi,H)}\|\beta\circ(\rho\otimes\pi')(p^n_G)\circ\alpha\|_{\mathcal{L}(V\otimes H)},$$
o{\`u} le supremum est pris parmi les repr{\'e}sentations $\pi$ unitaires de $\Gamma$ et
$\pi'$ est la repr{\'e}sentation de $G$ induite de $\pi$.\\
Donc,

                           $$\|p^n_{\Gamma}\|_{\mathcal{A}_{\Gamma}} \leq\sup\limits_{(\pi,H)}\mathrm{max}(\|\beta\circ\rho(p^n_G)\circ\alpha\|_{\mathrm{End}(V)}, \|\beta\circ(\rho\otimes\pi'_1)(p^n_G)\circ\alpha\|_{\mathcal{L}(V\otimes
                    H'_1)}),$$
o{\`u} on a {\'e}crit $\pi'=\pi'_0\oplus\pi'_1$ avec $\pi'_0$ {\'e}quivalente {\`a}
$1_G$ et $\pi'_1$ ne contenant pas $1_G$. Mais
$\lim_{n\rightarrow\infty}\beta\circ\rho(p^n_G)\circ\alpha=\mathrm{Id}_V$
et
$\lim_{n\rightarrow\infty}\beta\circ(\rho\otimes\pi'_1)(p^n_G)\circ\alpha=0$, donc
la suite $(p^n_{\Gamma})_n$ est une suite de Cauchy pour la norme
$\|  \|_{\mathcal{A}_{\Gamma}}$.

Soit $p_{\Gamma}$ sa limite dans $\in\mathcal{A}_{\Gamma}$.\\
On a alors l'{\'e}galit{\'e}:
$$(\rho\otimes\pi)(p_{\Gamma})=\beta\circ(\rho\otimes\pi')(p_{G})\circ\alpha.$$

\sloppy D'autre part, le fait d'avoir une fonction $f$ dans $C_c(G)$ qui v{\'e}rifie
$\sum_{\gamma\in\Gamma}f(g\gamma)\rho(\gamma)=\rho(g^{-1})$
implique que $\alpha(V\otimes H^{\Gamma})\subset V\otimes
(H')^G$ et que $\beta(V\otimes (H')^G)\subset
V\otimes H^{\Gamma}$. Comme de plus
$\beta\circ\alpha=\mathrm{Id}_{V\otimes H}$, alors $(\rho\otimes\pi)(p_{\gamma})$ est bien la projection de $V\otimes
H$ sur le sous-espace $V\otimes H^{\Gamma}$.\\

Montrons maintenant qu'il existe une fonction $f$ {\`a} support compact
sur $G$ et v{\'e}rifiant l'{\'e}quation (\ref{partition}).\\
Soit $p:G\rightarrow G/\Gamma$ la projection canonique. On doit
montrer qu'il existe une fonction $f$ continue {\`a} support compact telle
que, pour tout $x\in G/\Gamma$,
$$\sum\limits_{g\in p^{-1}(x)} f(g)\rho(g)=\mathrm{Id}_V.$$
Soit $(U_i)_{i=1,...,q}$ un
recouvrement ouvert de $G/\Gamma$ (que l'on peut prendre fini car
$G/\Gamma$ est compact) tel que $p^{-1}(U_i)$ soit hom{\'e}omorphe {\`a}
$U_i\times\Gamma$ et soient, pour tout $i=1,...,q$, $s_i:p^{-1}(U_i)\rightarrow
U_i\times\Gamma$ les hom{\'e}omorphismes correspondants. Comme $\rho$ est une repr{\'e}sentation irr{\'e}ductible de
$\Gamma$, l'ensemble des $\rho(\gamma)$ avec $\gamma$ parcourant
$\Gamma$, engendre $\mathrm{End}(V)$ (th{\'e}or{\`e}me de Burnside \cite[Chapter XVII, Corollary
3.3]{Lang}); si $m$ est la dimension de $V$,
on peut trouver un sous-ensemble $\Delta$ de $\Gamma$ de cardinal $m^2$ tel que
l'ensemble $\{\rho(\gamma) | \gamma\in\Delta\}$ forme une base de $\mathrm{End}(V)$. Il est clair alors que pour
tout $i=1,...,q$, il existe une fonction continue {\`a} support compact
$f_i:U_i\times\Gamma\rightarrow \mathbb{C}$ telle que
$$\sum\limits_{\gamma\in\Delta}f_i(u,\gamma)\rho(\gamma)=\rho(u^{-1}),$$
pour tout $u\in U_i$.\\
Pour tout $i=1,...,q$, soit $\tilde{f_i}=f_i\circ s_i$ et soit $(\delta_{i})_{i=1,...,q}$ une partition de l'unit{\'e} associ{\'e}e {\`a}
$(U_i)_{i=1,...,q}$.\\
On pose pour tout $g\in G$,
$$f(g)=\sum\limits_{i=1}^{q}(\delta_i\circ p)(g)\tilde{f_i}(g).$$
Comme $\tilde{f_i}$ appartient {\`a} $C_c(p^{-1}(U_i))$ pour tout
$i$, alors $(\delta_i\circ p)\tilde{f_i}$ est aussi à support compact sur $p^{-1}(U_i)$ et la fonction $f$
ainsi d{\'e}finie est {\`a} support compact sur $G$ et v{\'e}rifie
l'{\'e}quation (\ref{partition}).

\end{proof}

\section{Cas des groupes de Lie semi-simples}

Tout au long de cette section $G$ sera un groupe de Lie r{\'e}el connexe semi-simple de
centre fini. On ne consid{\`e}re que le cas o{\`u} $G$ est alg{\'e}brique et
simplement connexe (c'est-{\`a}-dire que tout rev{\^e}tement algébrique de
$G$ est isomorphe {\`a} $G$) de fa{\c c}on {\`a} ce que $G$ soit produit direct de
ses facteurs simples \cite[Proposition I.1.4.10]{Margulis}. On suppose en plus que tout facteur direct simple de $G$ est soit de rang au moins $2$, soit
localement isomorphe {\`a} $Sp(n,1)$ pour $n\geq 2$ ou {\`a} $F_{4(-20)}$. Le groupe $G$ a alors la propri{\'e}t{\'e} (T) de Kazhdan usuelle (c'est-{\`a}-dire que $1_G$ est isol{\'e}e dans
le dual unitaire de $G$) \cite{de la Harpe-Valette} et v{\'e}rifie, en
fait, une propri{\'e}t{\'e} plus forte de d{\'e}croissance uniforme de coefficients de matrice de repr{\'e}sentations unitaires \cite{Cowling}.\\ Soit $\rho:G\rightarrow
\mathrm{Aut(V)}$ une repr{\'e}sentation irr{\'e}ductible de dimension
finie de $G$ dans un espace vectoriel complexe $V$ de dimension
$m$. On consid{\`e}re le complexifi{\'e} $G(\mathbb{C})$ de $G$ et le complexifi{\'e}
de l'alg{\`e}bre de Lie de $G$, $\mathfrak{g}$, que l'on note
$\mathfrak{g}_\mathbb{C}$.  On notera de la m{\^e}me fa{\c c}on la repr{\'e}sentation
de $\mathfrak{g}$ d{\'e}finie par $\rho$ sur $V$ et $\rho$ elle m{\^e}me.
Soit $\mathfrak{u}$ une forme r{\'e}elle compacte de
$\mathfrak{g}_\mathbb{C}$ compatible avec $\mathfrak{g}$ et soit $U$ le
sous-groupe de Lie connexe de $G(\mathbb{C})$ qui a pour alg{\`e}bre de Lie
$\mathfrak{u}$. Le groupe $U$ est un sous-groupe compact maximal
de $G(\mathbb{C})$ qui est invariant par la conjugaison complexe sur
$G(\mathbb{C})$ \cite{Knapp2}. Soit $K$ le sous-groupe compact maximal de
$G$ donn{\'e} par $U\cap G$. On consid{\`e}re sur $V$ un produit hermitien
invariant par l'action de $U$ (unique {\`a} constante pr{\`e}s car $V$ est
irr{\'e}ductible), c'est-{\`a}-dire tel que $\rho_{\mathbb{C}}(U)$ soit contenu
dans les matrices unitaires de $\mathrm{GL}(V)$, o{\`u} $\rho_{\mathbb{C}}$
est le complexifi{\'e} de $\rho$. Pour un {\'e}l{\'e}ment $M\in
\mathrm{End}(V)$, notons $M^{*}$ son adjoint par rapport {\`a} ce
produit hermitien.

  On consid{\`e}re toujours la norme d'op{\'e}rateur
sur $\mathrm{End}(V)$ que l'on notera $\|.\|_ {\mathrm{End}(V)}$
et une mesure de Haar sur $G$, $dg$ pour $g\in G$. \\
Soit $\omega:G\rightarrow\mathrm{Aut}(V')$ une repr{\'e}sentation
fid{\`e}le de $G$ qui contient $\rho$ et telle que $V'$ soit muni
d'un produit hermitien de sorte que $\omega_{\mathbb{C}}(U)$ soit
contenu dans le groupe unitaire de $\mathrm{GL}(V')$. On d{\'e}finit
une longueur $l$ sur $G$ (c'est-{\`a}-dire une fonction sur $G$ {\`a}
valeurs dans $\mathbb{R}^{+}$ telle que
 $l(1)=0$, $l(gh)\leq l(g)+l(h)$ et $l(g)=l(g^{-1})$, $\forall g,h\in G$) de la fa{\c c}on suivante:
 $$l(g)=\log(\max(\|\omega(g)\|_{\mathrm{End}(V)},\|\omega(g^{-1})\|_{\mathrm{End}(V)})), \forall g\in G.$$
  \sloppy Cette longueur d{\'e}finit une semi-m{\'e}trique
$d$ sur $G$ donn{\'e}e par $d(g,x)=l(g^{-1}x)$, pour $g,x\in G$. Soit $B_q=\{g\in G | l(g)\leq q$\} pour tout $q\in\mathbb{N}$.\\

Le but de cette section est de d{\'e}montrer que le groupe $G$ a la
propri{\'e}t{\'e} (T) tordue par $\rho$.\\

Soient $\mathcal{A}_G$, $\mathcal{A'}_G$ et $\mathcal{A''}_G$
d{\'e}finies comme
dans la section pr{\'e}c{\'e}dente et soit $\Theta$ le morphisme d'alg{\`e}bres de Banach
$$\Theta:\mathcal{A}_G\rightarrow\mathcal{A'}_G\oplus\mathcal{A''}_G.$$
\begin{Theorem}\label{propo}
Le morphisme d'alg{\`e}bres de Banach $\Theta$ est un isomorphisme.
\end{Theorem}

\begin{Remark}\label{injectivit{\'e}}
Si on munit $\mathcal{A'}\oplus\mathcal{A''}$ de la norme donn{\'e}e
par: $\|(x,y)\|=\max(\|x\|_{\mathcal{A'}},\|y\|_{\mathcal{A''}})$
pour $x\in\mathcal{A'}$, $y\in\mathcal{A''}$, alors $\Theta$ est
un morphisme d'alg{\`e}bres de Banach isom{\'e}trique et donc pour prouver
le th{\'e}or{\`e}me il suffit de d{\'e}montrer que l'image est dense. En
effet, toute repr{\'e}sentation unitaire de $G$, $(\pi,H_{\pi})$, peut
s'{\'e}crire comme somme directe de deux sous-repr{\'e}sentations:
    la partie de $\pi$ qui ne contient pas de vecteurs invariants non nuls (et donc qui ne contient pas $1_G$)
    que l'on va noter $\pi_1$, et la partie de $\pi$ qui est multiple de $1_G$, not{\'e}e $\pi_0$. Alors, pour tout $f\in C_c(G)$,
    \begin{align*}
    \|f\|_{\mathcal{A}}&=\sup\limits_{\pi}(\max(\|(\rho\otimes\pi_1)(f)\|_{\mathcal{{L}}(V\otimes
                    H_{\pi_1})},\|\rho\otimes\pi_0 (f)\|_{\mathcal{L}(V\otimes H_{\pi_0})}))\\
               &=\mathrm{max}(\sup\limits_{\pi\nsupseteq 1_G}\|(\rho\otimes\pi)(f)\|_{\mathcal{{L}}(V\otimes
                    H)},\|\rho(f)\|_{\mathrm{End}(V)})\\
            &=\|\Theta(f)\|_{\mathcal{A'}\oplus\mathcal{A''}}.
            \end{align*}
            \end{Remark}

\begin{Remark}\label{Burnside}
On remarque aussi que $\rho$ d{\'e}finit une repr{\'e}sentation
irr{\'e}ductible et fid{\`e}le, car isom{\'e}trique, de $\mathcal{A''}$ dans
$\mathrm{End}(V)$. On a alors que $\mathcal{A''}=\mathrm{End}(V)$
(th{\'e}or{\`e}me de Burnside cf. \cite[Chapter XVII, Corollary
3.3]{Lang}).
\end{Remark}

\sloppy On a le lemme suivant:
    \begin{Lemma}\label{1}
    Il existe une matrice $E$ non nulle dans $\mathrm{End}(V)$ et il existe une suite de fonctions continues {\`a} support
    compact sur $G$, $(f_{n})_{n\in\mathbb{N}}$ telles que,
    $\lim\limits_{n\rightarrow\infty}\rho(f_{n})=E$ et $\lim\limits_{n\rightarrow\infty}\sup\limits_{\pi}\|(\rho\otimes\pi)(f_{n})\|_{\mathcal{L}(V\otimes H_\pi)}=0$,
    o{\`u} le supremum est pris parmi les repr{\'e}sentations unitaires de $G$ qui ne contiennent pas la repr{\'e}sentation triviale.
    \end{Lemma}
    Montrons d'abord que le lemme \ref{1} implique le Th{\'e}or{\`e}me
    \ref{propo}.\\
    Soit
    $(f_{n})_{n\in\mathbb{N}}$ la suite donn{\'e}e par le lemme \ref{1} et $E\in \mathrm{End}(V)$ non nulle telle que
    $\lim\limits_{n\rightarrow\infty}\rho(f_{n})=E$. La suite $(f_{n})_{n\in\mathbb{N}}$ converge dans $\mathcal{A}$ car
    $$\|f_n\|_{\mathcal{A}}=\mathrm{max}(\|\rho(f_n)\|_{\mathrm{End}(V)}, \sup\limits_{\pi\nsupseteq 1_G}\|(\rho\otimes\pi)(f_n)\|_{\mathcal{{L}}(V\otimes
                    H_\pi)}),$$
    et comme $\lim\limits_{n\rightarrow\infty}\sup\limits_{\pi\nsupseteq 1_{G}}\|(\rho\otimes\pi)(f_{n})\|_{\mathcal{L}(V\otimes H_\pi)}=0$,
    $(f_{n})_{n\in\mathbb{N}}$ est une suite de Cauchy pour $\|.\|_{\mathcal{A}}$ et donc elle converge dans $\mathcal{A}$.
    Soit $p$ la limite de $f_{n}$ quand $n$ tend vers l'infini. On a que $\rho(p)=E$ et pour toute repr{\'e}sentation unitaire
    $\pi$ de $G$, qui n'a pas de vecteurs invariants non nuls, $(\rho\otimes\pi)(p)=0$, donc $\Theta(p)=(0,E)$.
    Soit $\mathcal{S}$ l'id{\'e}al bilat{\`e}re engendr{\'e} par $p$ dans $\mathcal{A}$ et $Q$ l'id{\'e}al bilat{\`e}re engendr{\'e} par $E$
    dans $\mathrm{End}(V)$. Comme $\rho$ d{\'e}finit un morphisme surjectif de $\mathcal{A}$ vers $\mathrm{End}(V)$
    (cf. Remarque \ref{Burnside}), on a que $\Theta(\mathcal{S})$ contient $0\oplus Q$.
    Or,  tout id{\'e}al bilat{\`e}re non nul de $\mathrm{End}(V)$ est {\'e}gal {\`a} $\mathrm{End}(V)$
    tout entier. On a donc que $0\oplus\mathrm{End}(V)=0\oplus\mathcal{A''}$ est contenu dans
    $\Theta(\mathcal{S})$ et donc dans $\Theta(\mathcal{A})$. Ceci montre alors que
    $\Theta(\mathcal{A})=\mathcal{A'}\oplus\mathcal{A''}$, car $\mathcal{A}\rightarrow\mathcal{A'}$ est d'image dense,
    et donc que $\Theta$ est un isomorphisme en tenant compte de la remarque
    \ref{injectivit{\'e}}.

\begin{proof}
 
On va maintenant montrer le lemme \ref{1}. On remarque d'abord le lemme suivant:
\begin{Lemma}
    Soit $f$ une fonction {\`a} support compact sur $G$ et $(\pi,H_{\pi})$ une repr{\'e}sentation unitaire de $G$.
    On a l'in{\'e}galit{\'e} suivante:
      \begin{equation}\label{ineq}
        \|(\rho\otimes\pi)(f)\|_{\mathcal{{L}}(V\otimes
                    H_\pi)}\leq\sup
               m^2\int_{G}|f(g)|\|\rho(g)\|_{\mathrm{End}(V)}|\langle\pi(g)\xi,\eta\rangle|dg,
               \end{equation}
où le supremum est pris parmi les vecteurs unitaires $\xi, \eta$ de $H_{\pi}$.
        \end{Lemma}

\begin{proof} En effet,
\begin{align*}
 \|(\rho\otimes\pi)(f)\|_{\mathcal{L}(V\otimes H_{\pi})}&\leq\sup\limits_{
\begin{array}{c}
 \scriptstyle x, y\in V\otimes H_{\pi}\\
\scriptstyle \|x\|=\|y\|=1
\end{array}
}\int_G|f(g)||<(\rho\otimes\pi)(g)x,y>|dg.\\
\end{align*}
Or, pour tout $g\in G$ et pour tous vecteurs unitaires $x$ et $y$ de $V\otimes H_{\pi}$, il existe $\xi,\eta\in H_{\pi}$ tels que $\|\xi\|,\|\eta\|\leq m\|\rho(g)\|_{\mathrm{End}(V)}^{\frac{1}{2}}$ et
$$<(\rho(g)\otimes\pi(g))x,y>=<\pi(g)\xi,\eta>.$$
En effet, si on écrit $x=\sum\limits_{i}a_iv_i\otimes\alpha_i$ et $y=\sum\limits_{j}b_jv_j\otimes\alpha_j$, où $\{v_i\}_i$ et $\{\alpha_i\}_i$ sont des bases hilbertiennes de $V$ et de $H_{\pi}$ respectivement, alors
\begin{align*}
 <\rho(g)\otimes\pi(g)x,y>&=\sum\limits_{i,j}\overline{a_i}b_j<(\rho(g)\otimes\pi(g))(v_i\otimes\alpha_i),v_j\otimes\alpha_j>\\
&=\sum\limits_{i,j}\overline{a_i}b_j<\pi(g)\alpha_i,<\rho(g)v_i,v_j>\alpha_j>\\
&=<\pi(g)\xi,\eta>,
\end{align*}
où 
\begin{align*}
 \xi=\sum\limits_{i,j}a_i\overline{<\rho(g)v_i,v_j>}^{\frac{1}{2}}\alpha_i &\quad\hbox{et}\quad\xi=\sum\limits_{i,j}b_j<\rho(g)v_i,v_j>^{\frac{1}{2}}\alpha_j,
\end{align*}
ce qui implique 
\begin{align*}
\|\xi\|&\leq \sum\limits_{i}|a_i|\|\alpha_i\|\sum\limits_{j}\|\rho(g)\|_{\mathrm{End}(V)}^{\frac{1}{2}}\\
&\leq m\|\rho(g)\|_{\mathrm{End}(V)}^{\frac{1}{2}}
\end{align*}
 et de même $$\|\eta\|\leq m\|\rho(g)\|_{\mathrm{End}(V)}^{\frac{1}{2}}.$$\\
En particulier, on a que pour tout $g\in G$ et pour tous vecteurs unitaires $x$ et $y$ de $V\otimes H_{\pi}$, il existe $\xi',\eta'\in H_{\pi}$ tels que 
$$<(\rho(g)\otimes\pi(g))x,y>=m^2\|\rho(g)\|_{\mathrm{End}(V)}<\pi(g)\xi',\eta'>,$$
et $\|\xi'\|,\|\eta'\|\leq 1$ (il suffit de prendre $\xi'=\frac{\xi}{m\|\rho(g)\|^{\frac{1}{2}}}$ et $\eta'=\frac{\eta}{m\|\rho(g)\|^{\frac{1}{2}}}$).\\
Donc,
\begin{align*}
 \|(\rho\otimes\pi)(f)\|_{\mathcal{L}(V\otimes H_{\pi})}&\leq m^2\sup\limits_{
\begin{array}{c}
                            \scriptstyle\xi,\eta\in H_{\pi}\\
                            \scriptstyle\|\xi\|,\|\eta\|\leq 1
                            \end{array}
                            }
\int_{G}|f(g)|\|\rho(g)\|_{\mathrm{End}(V)}|\langle\pi(g)\xi,\eta\rangle|dg.
\end{align*}

\end{proof}

On veut utiliser la d{\'e}croissance uniforme des coefficients de
matrice des repr{\'e}sentations unitaires, ne contenant pas la
triviale, donn{\'e}e par le th{\'e}or{\`e}me suivant (voir \cite[Corollaire
2.4.3 et Th{\'e}or{\`e}me 2.5.3]{Cowling}, \cite[corollaire 2.7 et
proposition 6.3]{Howe}, \cite[proposition 2.7 et th{\'e}or{\`e}me
4.11]{Oh}):

\begin{Theorem}\label{Howe}
    Soit $G$ un groupe de Lie r{\'e}el semi-simple connexe {\`a} centre fini tel que, tout sous-groupe distingu{\'e}
    $G_{i}\neq 1$ soit tel que $\mathrm{rang}_{\mathbb{R}}(G_{i})\geq 2$, ou $G_i=Sp(n,1)$ ou $G_i=F_{4(-20)}$ et soit $K$ un sous-groupe compact maximal de $G$.
    Alors il existe une fonction continue $K$-bi-invariante $\phi$ sur $G$ {\`a} valeurs dans $\mathbb{R}^+$ qui tend vers z{\'e}ro
    {\`a} l'infini et telle que, pour toute repr{\'e}sentation unitaire $\pi$ de $G$ dans un espace de Hilbert $H_{\pi}$,
    qui ne contient pas de vecteurs invariants non nuls, et pour tous vecteurs unitaires $\xi$, $\eta$
    dans $H_{\pi}$, on a l'estimation suivante:
    $$\forall g\in G,
    |\langle\pi(g)\xi,\eta\rangle|\leqslant\phi(g)\delta_{K}(\xi)\delta_{K}(\eta)$$
o{\`u} $\delta_{K}(v)=(\mathrm{dim} \langle Kv\rangle)^{1/2}\in
\mathbb{N}\cup\{\infty\}$ et $\langle Kv\rangle$ est le
sous-espace de $V$ engendr{\'e} par l'action de $K$ sur $v$, pour
$v\in H_{\pi}$.
\end{Theorem}

On note $\widehat{K}$ l'ensemble des classes d'{\'e}quivalence de
repr{\'e}sentations irr{\'e}ductibles
    de $K$.\\
 On rappelle que toute  repr{\'e}sentation $(\mu,H_{\mu})$ de $K$ s'{\'e}crit
    comme somme directe de repr{\'e}sentations irr{\'e}ductibles. L'espace
    $H_{\mu}$ s'{\'e}crit alors comme une somme directe de la forme:
    $$H_{\mu}=\bigoplus_{[\sigma]\in\widehat{K}} H_{\sigma}^{\oplus
    r_{\sigma}},$$
   o{\`u} $[\sigma]$ est la classe de la repr{\'e}sentation $(\sigma,H_{\sigma})$ dans $\widehat{K}$
    et $r_{\sigma}$ est sa multiplicit{\'e} dans la d{\'e}composition de
    $\mu$.\\
    Le sous-espace $ H_{\sigma}^{\oplus
    r_{\sigma}}$ de $H_{\mu}$ est alors appel{\'e} la composante
    \emph{$\sigma$-typique} de $\mu$.\\
    Si $(\sigma,H_{\sigma})$ est une repr{\'e}sentation irr{\'e}ductible de dimension $n_{\sigma}$, la projection $P_{\sigma}:H_{\mu}\rightarrow H_{\mu}$
    sur la partie $\sigma$-typique de $\mu$ est donn{\'e}e par:
    \begin{equation}\label{projection}
    P_{\sigma}=n_{\sigma}\mu(\chi_{\sigma^*})\\
      \end{equation}
    o{\`u} $\chi_{\sigma}$ est le caract{\`e}re de $\sigma$ et $\chi_{\sigma^*}(t)=\overline{\chi_{\sigma}(t)}=\chi_{\sigma}(t^{-1})$
    est le caract{\`e}re de sa repr{\'e}sentation contragr{\'e}diente dans l'espace dual de
    $H_{\sigma}$ (cf. \cite[Chapitre 2, Partie I]{Serre}).\\

 Soit $\mathcal{I}\subset \widehat{K}$ l'ensemble des $K$-types de $V$, c'est-{\`a}-dire l'ensemble
    des repr{\'e}sentations irr{\'e}ductibles de $K$ qui apparaissent dans la d{\'e}composition de $(\rho\vert_K,V)$
    en somme directe de repr{\'e}sentations irr{\'e}ductibles. 
    \\

Pour toute repr{\'e}sentation $\pi$ de $G$, on note
$\mathcal{J}_{G}(\pi)$ l'ensemble des $K$-types de $\pi$ vue
comme repr{\'e}sentation de $G$. Avec cette notation $\mathcal{I}=\mathcal{J}_{G}(\rho)$.\\
En particulier, si $\pi$ est une repr{\'e}sentation de $G\times G$, on
notera $\mathcal{J}_{G\times G}(\pi)\subset \widehat{K}\times \widehat{K}$ l'ensemble
des $(K\times K)$-types de $\pi$ vue comme repr{\'e}sentation de $G\times G$.\\

    On consid{\`e}re la repr{\'e}sentation r{\'e}guli{\`e}re $L\times R$ de $G\times G$ sur $C_{c}(G)$, qui
    est donn{\'e}e par la formule:
    $$L\times R:G\times G \rightarrow \mathcal{L}(C_{c}(G)),\,\,(L\times R)(t,t')f(g)=f(t^{-1}gt').$$

\begin{Def}
Soit $f$ une fonction continue {\`a} support compact sur $G$.
L'ensemble des $K$-types {\`a} gauche et {\`a} droite de $f$ est
l'ensemble des classes de repr{\'e}sentations irr{\'e}ductibles de
$K\times K$ qui apparaissent dans la d{\'e}composition de $f\in
C_c(G)$ quand on d{\'e}compose la repr{\'e}sentation r{\'e}guli{\`e}re $L\times R$,
restreinte {\`a} $K\times K$, en somme directe de repr{\'e}sentations
irr{\'e}ductibles.
\end{Def}

\begin{Lemma}\label{K-types}
Il existe une fonction $\phi$ continue sur
$G$, bi-invariante par $K$, {\`a} valeurs dans $\mathbb{R}^+$, qui
tend vers z{\'e}ro {\`a} l'infini et telle que pour toute fonction $f$
continue {\`a} support compact sur $G$ ayant des $K$-types {\`a} gauche et
{\`a} droite contenus dans $\mathcal{J}_{G\times G}(\rho\otimes\rho^*)$ et pour toute repr{\'e}sentation
unitaire de $G$, $\pi$, sans vecteurs invariants non nuls, on a:

  $$\|(\rho\otimes\pi)(f)\|_{\mathcal{L}(V\otimes H_{\pi})}\leq\sup\limits_{\xi,\eta\in H_{\pi}}
        m^2\int_{G}|f(g)|\|\rho(g)\|_{\mathrm{End}(V)}\phi(g)\delta_{K}(\xi)\delta_{K}(\eta)dg,$$

o{\`u} le sup est pris parmi les vecteurs $\xi,\eta\in H_{\pi}$
unitaires ayant des $K$-types appartenant {\`a}
$\mathcal{J}_G(\rho\otimes\rho^*)$.
\end{Lemma}

\begin{proof}
On va d'abord montrer que le supremum dans l'inégalité (\ref{ineq}) peut être pris parmi les vecteurs unitaires $\xi, \eta\in H_{\pi}$ ayant des $K$-types appartenant à $\mathcal{J}_G(\rho\otimes\rho^*)$.\\
Soit $(\mu,H_{\mu})$ une repr{\'e}sentation unitaire de $G$ et $f$ une
fonction continue {\`a} support compact sur $G$. On remarque tout d'abord
que pour tout $\xi$,$\eta\in H_{\mu}$ et pour toute $\sigma\in\widehat{K}$, la projection du vecteur $\mu(f)\xi$
sur la composante $\sigma$-typique de $H_{\mu}$ est {\'e}gale {\`a}
$$P_{\sigma}(\mu (f)\xi)=n_{\sigma}\mu (\chi_{\sigma^*} * f)\xi,$$
et $\chi_{\sigma^*} * f$ est exactement la projection de $f$ sur la
composante $\sigma$-typique de la repr{\'e}sentation r{\'e}guli{\`e}re gauche $L$ de
$G$ sur $C_c(G)$. Donc $\chi_{\sigma^*} *f $
est non nul si et seulement si $\sigma$ est un
$K$-type {\`a} gauche de $f$. De m{\^e}me,
$$P_{\sigma}(\mu (f)^*\eta)=n_{\sigma}\mu(f*\chi_{\sigma})^*\xi,$$
et $f*\chi_{\sigma}$ {\'e}tant la projection de $f$ sur la composante
$\sigma^*$-typique de la repr{\'e}sentation r{\'e}guli{\`e}re droite $R$ de $G$ sur $C_c(G)$, il est non nul si et
seulement si $\sigma^*$ est un $K$-type {\`a} droite de $f$.\\

 Soient
$\mathcal{K}_{f,L}$ et $\mathcal{K}_{f,R^*}$ les deux
sous-ensembles de $\widehat{K}$ d{\'e}finis par:
\begin{align*}
&\mathcal{K}_{f,L}=\{\sigma\in\widehat{K}| \,\sigma  \hbox{ est $K$-type {\`a} gauche de $f$}\}\\
&\mathcal{K}_{f,R^*}=\{\sigma\in\widehat{K}|\, \sigma^* \hbox{est
$K$-type {\`a} droite de $f$}\}.
\end{align*}

Comme les projections $P_{\sigma}$ sur les espaces $\sigma$-typiques sont
des projections orthogonales (Lemme de Schur \cite{Serre}), alors:
           $$\|\mu(f)\|_{\mathcal{{L}}(H_{\mu})}=\sup\limits_{z,y}|\langle\mu(f)z,y\rangle|,$$
            o{\`u} $z$ et $y$ parcourent les vecteurs unitaires de $H_{\mu}$ tels que l'ensemble des $K$-types de $z$ est contenu
        dans $\mathcal{K}_{f,R^*}$ et l'ensemble des $K$-types de $y$ est contenu dans $\mathcal{K}_{f,L}$.\\

Soit $(\pi,
H)$ une repr{\'e}sentation unitaire de $G$. On consid{\`e}re le produit tensoriel $(\rho\otimes\pi,V\otimes H)$.\\
Supposons maintenant que  les $K$-types {\`a} gauche de $f$ soient
contenus dans l'ensemble des $K$-types de $V$ et les $K$-types {\`a}
droite de $f$ soient contenus dans l'ensemble des $K$-types de
$V^*$. Alors, $\mathcal{K}_{f,L}\subset\mathcal{I}$ et comme
             $$\mathcal{K}_{f,R^*}\subset\{\sigma\in\widehat{K}|\, \sigma^* \hbox{est
$K$-type de $V^*$}\},$$
on a aussi que $\mathcal{K}_{f,R^*}\subset\mathcal{I}$. D'o{\`u} l'in{\'e}galit{\'e},
$$\|(\rho\otimes\pi)(f)\|_{\mathcal{{L}}(V\otimes
H)}\leq\sup\limits_{z,y}
|\int_{G}f(g)\langle(\rho\otimes\pi)(g)z,y\rangle dg|,$$ o{\`u} le
supremum est pris parmi les vecteurs unitaires $z$, $y\in V\otimes
H$ ayant des $K$-types appartenant {\`a} l'ensemble des $K$-types de
$V$.\\

De plus, si $H'$ est le sous-espace vectoriel de $H$ form{\'e} des
vecteurs dont les $K$-types sont parmi les $K$-types de $V\otimes
V^*$, tout vecteur de $H\otimes V$ dont les $K$-types sont parmi
ceux de $V$ appartient {\`a} $H'\otimes V$ (car
$\mathrm{Hom}_K(V\otimes V^*,H)=\mathrm{Hom}_K(V,V\otimes H)$, o{\`u}
$\mathrm{Hom}_K$ d{\'e}signe l'espace des morphismes de repr{\'e}sentations
$K$-invariants). Donc l'inégalité (\ref{ineq}) devient:
        $$\|(\rho\otimes\pi)(f)\|_{\mathcal{{L}}(V\otimes
        H)}\leq\sup\limits_{\xi,\eta}
                 m^2\int_{G}|f(g)|\|\rho(g)\|_{\mathrm{End}(V)}|\langle\pi(g)\xi,\eta\rangle|dg,$$
        o{\`u} $\xi$ et $\eta$ sont des vecteurs unitaires de $H$ qui ont des $K$-types appartenant {\`a}
    l'ensemble des $K$-types de $V\otimes V^*$, $\mathcal{J}_G(\rho\otimes\rho^*)$.\\

Consid{\'e}rons une fonction $\phi$ continue, $K$-bi-invariante sur
$G$ et {\`a} valeurs dans $\mathbb{R}^+$, qui tend vers z{\'e}ro
    {\`a} l'infini, donn{\'e}e par le th{\'e}or{\`e}me \ref{Howe}. Alors on a
  pour toute repr{\'e}sentation
unitaire $\pi$ qui ne contient pas la triviale:
 $$\|(\rho\otimes\pi)(f)\|_{\mathcal{{L}}(V\otimes H)}\leq
 m^2\sup\limits_{\xi,
 \eta}\int_{G}|f(g)|
 \|\rho(g)\|_{\mathrm{End}(V)}\phi(g)\delta_{K}(\xi)\delta_{K}(\eta)dg,$$
 o{\`u} $\xi$ et $\eta$ parcourent les vecteurs unitaires de $H$ qui sont $K$-finis et qui ont des $K$-types
 contenus dans $\mathcal{J}_G(\rho\otimes\rho^*)$.\\
\end{proof}

On cherche maintenant une
suite de fonctions $f_{n}\in C_{c}(G)$ ayant des $K$-types {\`a}
gauche et {\`a} droite contenus dans $\mathcal{J}_{G\times
  G}(\rho\otimes\rho^*)$. Pour simplifier la notation, on note $\mathcal{J}:=\mathcal{J}_{G\times
  G}(\rho\otimes\rho^*)$. On a le lemme
suivant:
    \begin{Lemma}\label{2}
    Il existe une matrice  non nulle  $E\in \mathrm{End}(V)$ et il existe une suite de fonctions
    $f_{n}\in C_{c}(G)$ ayant des $K$-types {\`a} gauche et {\`a} droite contenus dans $\mathcal{J}$
    et une constante positive $D$ telles que, pour tout entier $n$, le support de $f_{n}$ soit contenu dans
    $G\backslash B_{n}$, $\lim\limits_{n\rightarrow\infty}\rho(f_{n})=E$ et
    $$\int_{G} |f_{n}(g)|||\rho(g)||_{\mathrm{End}(V)}dg\leq
    D.$$\\
    \end{Lemma}

    \begin{proof}
On consid{\`e}re la d{\'e}composition de Cartan de $\mathfrak{g}$ donn{\'e}e
par la forme r{\'e}elle compacte $\mathfrak{u}$. $\mathfrak{g}$
s'{\'e}crit alors
$$\mathfrak{g}=\mathfrak{k}\oplus\mathfrak{p},$$
o{\`u} $\mathfrak{k}=\mathfrak{g}\cap\mathfrak{u}$ et
$\mathfrak{p}=\mathfrak{g}\cap i\mathfrak{u}$, et $\mathfrak{p}$
est non nul car $G$ n'est pas compact. Alors pour tout
$x\in\mathfrak{u}$, $\rho(x)$ est une matrice anti-hermitienne et
pour tout $x\in i\mathfrak{u}$, $\rho(x)$ est hermitienne.\\
Soit $X\in\mathfrak{p}$ non nul et $a:=\exp(X)$. Par cons{\'e}quent
$\rho(a)=\exp(\rho(X))$ est une matrice hermitienne (donc
diagonalisable dans une base orthonormale de $V$) {\`a} valeurs
propres strictement positives que l'on notera, sans tenir compte
des multiplicit{\'e}s,
$\nu_1,....,\nu_m$, o{\`u} $\nu_i\in\mathbb{R}^*_+$ pour tout $1\leq i\leq m$.\\
De plus, comme $\omega$ est une repr{\'e}sentation fid{\`e}le de $G$ qui
envoie $U$ dans les matrices unitaires de $\mathrm{End}(V')$,
$\omega(X)$ est une matrice hermitienne non nulle. On a que
$l(\exp(tX))=t\|\omega(X)\|_{\mathrm{End}(V')}$, car
$\|\omega(\exp(tX))\|_{\mathrm{End}(V')}=\exp(t\|\omega(X)\|_{\mathrm{End}(V')})$,
pour tout r{\'e}el $t$. Quitte {\`a} remplacer $a$ par $a^k$ pour $k$ un
entier assez grand,
on peut m{\^e}me supposer $l(a)\geq 2$.\\
Posons $a_n=a^n$, pour tout entier positif $n$. On a alors que
$l(a_n)=n\|\omega(X)\|_{\mathrm{End}(V')}=nl(a)$, donc
$l(a_n)\geq 2n$ et
$a_n$ appartient {\`a} $G\backslash B_n$.\\
De plus on a, pour tout $n$, $\rho(a_n)=\rho(a)^n$ est
diagonalisable sur $\mathbb{R}$ dans la m{\^e}me base que $\rho(a)$ et
ses valeurs propres sont ${\nu_1}^n,....,{\nu_m}^n$. Comme
$\|\rho(a_n)\|_{\mathrm{End}(V)}=\max\limits_{1\leq i\leq
m}(\nu_i^n)$ et que $\nu_i> 0$, pour tout $i=1...m$, on a
que
  $\frac{\rho(a_n)}{\|\rho(a_n)\|_{\mathrm{End}(V)}}$ tend vers une matrice
de $\mathrm{End}(V)$ non nulle, que l'on note $E'$.\\

Soit maintenant $f$ une fonction continue positive {\`a} support
compact sur $G$ telle que $\int_Gf(g)dg=1$ et telle que le support
de $f$ soit contenu dans
$B_1\cap \{g\in G |\,\, \|\rho(g)-\mathrm{Id}\|_{\mathrm{End}(V)}\leq \frac{1}{2}\}$.\\
Soit $f_n$ dans $C_c(G)$ d{\'e}finie de la fa{\c c}on suivante: pour tout
$g$ dans $G$
$$f_n(g)=\frac{f(a_n^{-1}g)}{\|\rho(a_n)\|_{\mathrm{End}(V)}}.$$
On a donc que $\mathrm{supp} f_n\subset a_n(\mathrm{supp} f)$ et
$f_n$ est dans $C_c(G)$, pour tout $n$. De plus, on a que pour
tout $g$ appartenant au support de $f_n$,
   $$l(g)\geq l(a_n)-1,$$
   donc le support de $f_n$ est contenu dans $G\backslash B_n$.\\
D'autre part, on voit facilement que $$\rho(f_n)=\frac{\rho(a_n)}{\|\rho(a_n)\|_{\mathrm{End}(V)}}\int_Gf(g)\rho(g)dg.$$\\
Posons $J:=\int_Gf(g)\rho(g)dg$. $J$ est une matrice inversible de
$\mathrm{End}(V)$. En effet, on a que
\begin{align*}
\|\int_Gf(g)\rho(g)dg-\mathrm{Id}_{\mathrm{End}(V)}\|_{\mathrm{End}(V)}&\leq \int_Gf(g)\|\rho(g)-\mathrm{Id}_{\mathrm{End}(V)}\|_{\mathrm{End}(V)}dg\\
&\leq \sup_{g\in \mathrm{supp}(f)}\|\rho(g)-\mathrm{Id}_{\mathrm{End}(V)}\|_{\mathrm{End}(V)}\\
&\leq\frac{1}{2}.
\end{align*}
 Par cons{\'e}quent, $\rho(f_n)$ tend vers $E'J$ qui
est encore une matrice non nulle de $\mathrm{End}(V)$.\\
 Par
ailleurs, on a que,
\begin{align*}
\int_Gf_n(g)\|\rho(g)\|&_{\mathrm{End}(V)}dg\\
&\leq\frac{1}{\|\rho(a_n)\|_{\mathrm{End}(V)}}\int_Gf(g)\|\rho(a_n)\|_{\mathrm{End}(V)}\|\rho(g)\|_{\mathrm{End}(V)}dg\\
                         &\leq\int_Gf(g)\|\rho(g)\|_{\mathrm{End}(V)}dg\\
                         &\leq\frac{3}{2}.
\end{align*}
On a donc trouv{\'e} une suite de fonctions $f_n$ dans $C_c(G)$, une
matrice $E=E'J$ dans $\mathrm{End}(V)$ non nulle et une contante
$D$ tels que $\lim\limits_{n\rightarrow\infty} \rho(f_n)=E$,
$\int_{G}|f_{n}(g)|\|\rho(g)\|_{\mathrm{End}(V)}dg\leq D$ et tels
que le support de $f_n$, pour tout $n$, soit contenu dans
$G\backslash B_n$. On va montrer qu'on peut
prendre les fonctions $f_n$ ayant des $K$-types, {\`a} droite et {\`a} gauche, contenus dans $\mathcal{J}$.\\

 Soit une fonction $f\in C_c(G)$. D'apr{\`e}s la formule (\ref{projection}),
    qui donne la projection sur les
    composantes $\sigma$-typiques
    d'une repr{\'e}sentation de $G$, pour $\sigma\in\widehat{K}$, la
    fonction $\tilde{f}$ d{\'e}finie par la formule suivante:
    \begin{align*}
\tilde{f}&=\sum\limits_{\phi_{1},\phi_{2}\in I}n_{\scriptstyle\phi_{1}\otimes\phi_{2}^*}(L\times R)( \chi_{\scriptstyle\phi_{1}^*\otimes\phi_{2}})\\
&=\sum\limits_{\phi_{1},\phi_{2}\in I}n_{\scriptstyle\phi_{1}}\chi_{\scriptstyle\phi_{1}^*} * f * n_{\scriptstyle\phi_{2}}\chi_{\scriptstyle\phi_{2}^*}.
    \end{align*}
 est la projection de $f$ sur les $K\times K$-types de $V\otimes V^*$.\\

On utilise ceci pour obtenir, pour tout $n$, une fonction
$\tilde{f_{n}}$ qui est donn{\'e}e par la projection de $f_{n}$ sur
les $K\times K$-types de $V\otimes V^*$. On a alors une suite de
fonctions $(\tilde{f_n})_{n\in\mathbb{N}}$ ayant des $K$-types {\`a}
gauche et {\`a} droite appartenant {\`a} $\mathcal{J}$.

On va maintenant v{\'e}rifier que la nouvelle suite satisfait les conditions du lemme \ref{2}.\\
L'application $\rho:C_{c}(G)\rightarrow \mathrm{End}(V)\backsimeq
V\otimes V^*$ est un morphisme de repr{\'e}sentations de $G\times G$
et, pour tous $t, t'$ dans $G$ le diagramme suivant commute:
    $$\xymatrix{
    C_c(G)\ar[d]_{(L\times R)(t,t')}\ar[r]^{\rho} & V\otimes V^*\ar[d]^{(\rho\otimes\rho^*)(t,t')}\\
    C_c(G)\ar[r]_{\rho} & V\otimes V^*\\}$$

Par cons{\'e}quent, comme
$\lim\limits_{n\rightarrow\infty}\rho(f_{n})=E$,
$\lim\limits_{n\rightarrow\infty}\rho(\tilde{f_{n}})$ est {\'e}gal {\`a}
la projection de $E$ sur les $(K\times K)$-types de $V\otimes V^*$,
et cette projection n'est rien d'autre que $E$
elle m{\^e}me, c'est {\`a} dire que $\lim\limits_{n\rightarrow\infty}\rho(\tilde{f_{n}})=E$.\\
De plus, on a que,
    \begin{align*}
        \int_{G}|\tilde{f_{n}}(g)|&
        \|\rho(g)\|_{\mathrm{End}(V)}dg\\
    &\leq\int_{G}\sum\limits_{\phi_{1},\phi_{2}\in I}n_{\phi_{1}}.n_{\phi_{2}^*}\int_{K\times K}
                            |\chi_{\phi_{1}^*}(t)\chi_{\phi_{2}}(t')||f_{n}(t^{-1}gt')|\|\rho(g)\|dtdt'dg\\
                        &\leq\sum\limits_{\phi_{1},\phi_{2}\in I}n_{\phi_{1}}.n_{\phi_{2}^*}\int_{K\times K}|
                            \chi_{\phi_{1}^*}(t)\chi_{\phi_{2}}(t')|\int_{G}|f_{n}(t^{-1}gt')|\|\rho(g)\|dgdtdt'\\
                        &\leq D\sum\limits_{\phi_{1},\phi_{2}\in I}n_{\phi_{1}}.n_{\phi_{2}^*}
                        \int_{K\times
                        K}|\chi_{\phi_{1}^*}(t)\chi_{\phi_{2}}(t')|dtdt'\\
                        &\leq D',
                        \end{align*}
     o{\`u} $D'$ est une constante qui ne d{\'e}pend pas de $n$.
     Comme le support de $\tilde{f_{n}}$ est contenu dans $K(\mathrm{supp} f_n)K$ pour tout $n$, et les $B_n$ sont invariants par l'action
     {\`a} gauche et {\`a} droite de $K$, le support de $\tilde{f_{n}}$ est contenu dans $G\backslash B_n$.
    \end{proof}

\sloppy Soient maintenant $E\in \mathrm{End}(V)$ non nulle et
$f_{n}\in C_{c}(G)$ tels que
$\lim\limits_{n\rightarrow\infty}\rho(f_{n})=E$, avec le support
de $f_n$ contenu dans
        $G\backslash B_n$, $f_n$ ayant des $K$-types {\`a} gauche et {\`a} droite contenus $\mathcal{J}$,
    et tels que $\int_{G}|f_{n}(g)|
        \|\rho(g)\|_{\mathrm{End}(V)}dg\leq D$, pour une constante $D$ qui ne d{\'e}pend que de $E$.\\
D'apr{\`e}s le lemme \ref{K-types}, si $\pi$ est une repr{\'e}sentation
unitaire de $G$ qui ne contient pas la triviale, on a que, pour
tout entier naturel $n$,
 $$\|(\rho\otimes\pi)(f_{n})\|_{\mathcal{{L}}(V\otimes H)}\leq
 m^2\sup\limits_{\xi,
 \eta}\int_{G}|f_{n}(g)|
 \|\rho(g)\|_{\mathrm{End}(V)}\phi(g)\delta_{K}(\xi)\delta_{K}(\eta)dg,$$
 o{\`u} $\xi$ et $\eta$ parcourent les vecteurs unitaires de $H$ qui ont des $K$-types
 contenus dans l'ensemble des $K$-types de $V\otimes V^*$ et $\phi$ est une fonction continue et positive sur $G$ qui s'annule {\`a} l'infini
 et qui ne d{\'e}pend ni de $f_n$ ni de $\pi$.
 \\
De plus, on a une constante positive $D$ telle que
$\int_{G}|f_{n}(g)|\|\rho(g)\|_{\mathrm{End}(V)}\leq D$ et, comme
le support de $f_n$ est contenu dans $G\backslash B_n$,
l'in{\'e}galit{\'e} au-dessus s'{\'e}crit:
\begin{equation}\label{ineq2}
    \|(\rho\otimes\pi)(f_{n})\|_{\mathcal{L}(V\otimes H)}\leq
    m^{2}D\sup\limits_{\xi,\eta}\sup\limits_{g\in G\backslash
    B_{n}} \phi(g)\delta_{K}(\xi)\delta_{K}(\eta),
 \end{equation}
o{\`u} $\xi$ et $\eta$ parcourent les vecteurs unitaires de $H$ qui
ont des $K$-types appartenant {\`a} l'ensemble des $K$-types
de $V\otimes V^*$, ensemble qui ne d{\'e}pend pas de $n$.\\
On veut montrer que le membre de droite de cette in{\'e}galit{\'e} tend
vers z{\'e}ro quand $n$ tend vers l'infini. Pour ceci, on a besoin du
lemme suivant qui assure que, pour tout $v\in H$ ayant des
$K$-types contenus dans un ensemble fix{\'e} $S$, la dimension du
sous-espace de $H$ engendr{\'e} par l'action de $K$ sur $v$, que l'on
note $\delta_{K}(v)$, est born{\'e}e en fonction de $S$.
 \begin{Lemma}\label{3}
 Soit $v\in W$, o{\`u} $(\mu,W)$ est une repr{\'e}sentation de $K$. Alors,
 $$\delta_{K}(v)=\mathrm{dim}\langle Kv \rangle\leq \sum\limits_{[\sigma]} (\mathrm{dim}\,\sigma)^2,$$
 o{\`u} la somme est prise parmi les $[\sigma]\in \widehat{K}$ qui sont des $K$-types de $v$.
 \end{Lemma}
 \begin{proof}
 Soit $C_r^*(K)$ la C*-alg{\`e}bre r{\'e}duite de $K$. Toute repr{\'e}sentation irr{\'e}ductible $\sigma$ de $K$, apparait
$\mathrm{dim}(\sigma)$ fois dans la d{\'e}composition de la
repr{\'e}sentation r{\'e}guli{\`e}re de $K$ en somme directe de
repr{\'e}sentations irr{\'e}ductibles (cf. \cite{Serre}). En fait,
l'application
\begin{align*}
C_r^*(K)&\stackrel{\simeq}{\rightarrow} \bigoplus\limits_{[\sigma]\in\widehat{K}}\mathrm{End}(H_{\sigma})\\
f&\mapsto (\sigma(f))_\sigma,
\end{align*}
est un isomorphisme de C*-alg{\`e}bres. Soit $\psi$ le morphisme de
$C_r^*(K)$ vers $W$ qui envoie $f$ dans $\mu(f)v\in \langle Kv
\rangle$. On a donc que,
\begin{align*}
\langle Kv \rangle&=\psi(C_r^*(K))\\
&=\psi(\bigoplus\limits_{\scriptstyle{[\sigma]\in \widehat{K},}\
                      }\mathrm{End}(H_\sigma)),
\end{align*}
o{\`u} la somme directe est prise parmi les $[\sigma]$
                      qui sont des $K$-types de $v$.\\
 On a alors que $$\mathrm{dim}\langle Kv \rangle\leq \sum\limits_{[\sigma]} (\mathrm{dim}\,\sigma)^2,$$
 o{\`u} la somme est prise parmi les $[\sigma]\in \widehat{K}$ qui sont des $K$-types de $v$.
  \end{proof}

Maintenant on est pr{\^e}t pour conclure. Le membre de droite de
(\ref{ineq2}) tend vers z{\'e}ro quand $n$ tend vers l'infini, car la
fonction $\phi$ tend vers z{\'e}ro {\`a} l'infini, et donc la norme de
$(\rho\otimes\pi)(f_{n})$
tend vers z{\'e}ro quand $n$ tend vers l'infini.\\
On a donc trouv{\'e} une matrice $E$ non nulle dans $\mathrm{End}(V)$
et une suite de fonctions $f_{n}$ dans $C_{c}(G)$ tels que
$\lim\limits_{n\rightarrow\infty}\rho(f_{n})=E$ et le supremum sur
toutes les repr{\'e}sentations unitaires $\pi$ de $G$, qui ne
contiennent pas la triviale, de
$\|(\rho\otimes\pi)(f_{n})\|_{\mathcal{L}(V\otimes H_\pi)}$ tend
vers z{\'e}ro quand $n$ tend vers l'infini. Ceci termine la
d{\'e}monstration du lemme \ref{1}, donc celle du th{\'e}or{\`e}me
\ref{propo}. Par la proposition \ref{evidente}, $G$ a alors la propriété (T) tordue par $\rho$, ce qui termine la démonstration du théorème \ref{intro} énoncé dans l'introduction.
\end{proof}


\begin{thebibliography}{99}
\bibitem[AW]{Akemann-Walter} C.A Akemann, M.E. Walter, \emph{Unbounded Negative Definite Functions}, Can. J. Math., Vol. XXXIII, No. 4, 1981, pp.862-871
    \bibitem[Di]{Dixmier} J. Dixmier, \emph{C*-algebras}, North-Holland Publ., Amsterdam 1977
        \bibitem[DK]{Duistermaat-Kolk} J. J. Duistermaat, J. A.C. Kok, \emph{Lie Groups}, Springer-Verlag, 1942
        \bibitem[Co]{Cowling} M. Cowling, \emph{Sur les coefficients des repr{\'e}sentations unitaires des groupes de Lis semi-simples, in:
    P. Eymard, J. Faraut, G. Schiffmann, and R. Takahashi, eds.}, Analyse Harmonique sur les Groupes de Lie II (S{\'e}minaire
    Nancy-Strasbourg 1976-78), Lecture Notes in Mathematics 739, Springer-Verlag, New York, 132-178, 1979
    \bibitem[FD]{Fell} J. M. G. Fell, R. S. Doran,
    \emph{Representations of *-Algebras, Locally Compact Groups, and
    Banach *-Algebraic Bundles- Volume 1- Basic Representation Theory
    of Groups and Algebras}, Acamedic Press, INC., 1988
    \bibitem[FH]{Fisher} D. Fisher, T. Hitchman, \emph{Strengthening
    Kazhdan's Property (T) by Bochner Methods}, preprint, arxiv:math.DG/0609663  
 \bibitem[Go]{Gomez} M.P. Gomez-Aparicio, \emph{Property (T) and tensor products by irreducible finite dimensional representations for $SL_n(\RR)$, $n\geq 3$}, à paraître dans AMS, Contemporary Mathematics
\bibitem[Ho]{Howe} R.E. Howe, \emph{On a notion of rank for unitary representations of the classical groups}, Harmonic    analysis and group representations, C.I.M.E, II Circlo, Palazzone della Scuola Normale Superiore, 223-331, 1980
    \bibitem[HT]{Howe-Tan} R.E. Howe, E.C Tan, \emph{Non-abelian Harmonic Analysis: applications to $SL(2,\mathbb{R})$},
    Springer-Verlag, 1992
    \bibitem[HV]{de la Harpe-Valette} P. de la Harpe, A. Valette, \emph{La propri{\'e}t{\'e} (T) de Kazhdan
    pour les groupes localement compacts}, Ast{\'e}risque 175, 1989
    \bibitem[Ka]{Kazhdan} D.A. Kazhdan, \emph{Connection of the dual space of a group with the structure of its closed subgroups},
    Functional Analysis and its Applications 1, 1967, 63-65
    \bibitem[Kn]{Knapp} A.W. Knapp, \emph{Representation theory of semi-simple groups}, Princeton University, 1986
    \bibitem[Kn2]{Knapp2} A.W. Knapp, \emph{Lie Groups Beyond an introduction}, Birkh{\"a}user, 1996
    
    \bibitem[La]{Lang} S. Lang, \emph{Algebra},
    Springer-Verlag, 2002
    \bibitem[Ma]{Margulis} G. A. Margulis, \emph{Discrete subgroups of
    semisimple Lie groups}, Ergebnisse der Mathematik und ihrer
    Grenzgebiete (3) [Results in Mathematics and Related Areas (3)],
    vol. 17 Springer-Verlag, Berlin, 1991
  \bibitem[Oh]{Oh} H. Oh, \emph{Uniform pointwise bounds for matrix coefficients of unitary representations and applications
    to Kazhdan constants,} Duke Math. J. 113, 133-192, 2002
    \bibitem[Se]{Serre} J.P. Serre, \emph{Repr{\'e}sentations lin{\'e}aires des groupes finis}, Paris: Hermann, 1971
    \bibitem[Va]{Valette} A. Valette, \emph{Minimal projections, integrable representations and property (T)},
    Arch. Math., Vol.43, 397-406, 1984


\end{thebibliography}
\end{document}